\documentclass{aptpub}
\usepackage{setspace}

\authornames{S. SAGITOV AND M.C. SERRA } 
\shorttitle{Skeletons of near-critical Bienaym\'{e}-Galton-Watson branching processes} 

\newcommand{\be}{\begin{equation}} \newcommand{\ee}{\end{equation}}
\usepackage{graphicx}
\numberwithin{equation}{section}  
\numberwithin{definition}{section}

\begin{document}

\title{Skeletons of near-critical Bienaym\'{e}-Galton-Watson branching processes} 

\authorone[Chalmers University of Technology and Gothenburg University]{Serik Sagitov}  
\authortwo[Center of Mathematics, Minho University]{Maria Concei\c{c}\~{a}o Serra} 

\begin{abstract}
Skeletons of branching processes are defined as trees of lineages characterized by an appropriate signature of future reproduction success. In the supercritical case a natural choice is to look for the lineages that survive forever \cite{O}. In the critical case it was earlier suggested \cite{S} to distinguish the particles with the total number of descendants exceeding a certain threshold. These two definitions lead to asymptotic representations of the skeletons as either pure birth process (in the slightly supercritical case) or critical birth-death processes (in the critical case conditioned on the total number of particles exceeding a high threshold value). The limit skeletons reveal typical survival scenarios for the underlying branching processes.

In this paper we consider near-critical Bienaym\'{e}-Galton-Watson processes and define their skeletons using marking of particles. If marking is rare, such skeletons are approximated by birth and death processes which can be subcritical, critical or supercritical. We obtain the limit skeleton for a sequential mutation model \cite{SS} and compute the density distribution function for the time to escape from extinction.
\end{abstract}

\keywords{Bienaym\'e-Galton-Watson process; decomposable multi-type branching process; birth and death process; escape from extinction} 

\ams{60J80}{92D25} 

\addressone{Mathematical Sciences, Chalmers and Gothenburg University, SE-41296 Gothenburg, Sweden. Email address: serik@chalmers.se} 
\addresstwo{Center of Mathematics, Minho University, Campus de Gualtar, 4710 - 057 Braga, Portugal. Email address: mcserra@math.uminho.pt} 

\section{Introduction} 

Imagine a population of viruses trying to establish itself in a new environment. Suppose the currently dominating type is nearly critical, in that its mean offspring number is close to one. One can think of two main factors which may lead to survival of this population: reproductive success or an advantageous mutation (a mutation producing new type of particles forming a strictly supercritical process). While a reproductive success is possible in the slightly supercritical case, `survival due to an advantageous mutation' is the only way to escape extinction for a slightly subcritical  branching system.

The typical survival scenarios of such branching processes can be studied in terms of the so-called skeleton trees formed by lineages characterized by an appropriate signature of future reproduction success. In the supercritical case a natural choice is to look for the lineages that survive forever \cite{O}. In the critical case it was earlier suggested \cite{S} to distinguish the particles with the total number of descendants exceeding a certain threshold. These two definitions lead to asymptotic representations of the skeletons as either pure birth process (in the slightly supercritical case) or critical birth-death processes (in the critical case conditioned on the total number of particles exceeding a high threshold value).

In this paper we suggest an alternative approach of defining a skeleton that relies on a random marking of the lineages in the family tree of a Bienaym\'{e}-Galton-Watson (BGW) process. The skeleton is then defined as the subtree formed by the infinite lineages together with the marked lineages.
In Section \ref{Ssc} we describe the skeleton of infinite lineages and recall the result from \cite{O} concerning a sequence of single type slightly-supercritical BGW process. It says that, conditioned on the event that the skeleton is not empty, the skeleton is approximated by the standard Yule process (a linear pure birth process).

In Section \ref{Scr} we consider an exactly critical BWG process with marking: each particle in the family tree is marked, independently of the others, with a small probability. Here the skeleton is formed by the lineages leading to a marked particle. We show that, conditioned on the event that it is not empty, the skeleton is approximated by a critical linear birth-death process.

In Section \ref{Snc} we state the main result of the paper concerning a sequence of near-critical BGW process with marked particles. The definition of skeleton is adapted in order to include both infinite lineages and marked lineages. The marking is also done in a more general way than in Section \ref{Scr}. Our main result states that, conditioned on the event that it is not empty, the skeleton is approximated by a linear birth-death process which can be either supercritical, critical or subcritical, depending on the parameters of the model. The proof of this main result, Theorem \ref{Thm1}, is given in Section 8.

In Section \ref{Sbm} a decomposable two-type BGW process with irreversible mutations, starting from a single wild type individual, is studied. Each daughter of a wild type individual becomes a mutant, independently of the others, with a small probability. We look at this Binomial mutation model as a particular case of the processes treated in Section \ref{Snc} by considering the following marking procedure: a wild type individual is marked if gives birth to at least one mutant daughter.

Section \ref{Ssm} deals with a sequential mutation model, considered in \cite{SS}, for a viral population that escapes extinction due to a sequence of irreversible mutations that lead to a target type. It is assumed that mutations appear according to a Binomial mutation model and results from Section \ref{Sbm} are used to obtain the asymptotic shape of the limit skeleton. Finally,  in Section 7 we address the important question of the time to escape from extinction in a sequential mutation model. By 'the time to escape from extinction' we mean the first generation where a particle of the target type appears. Due to the shape of the skeleton, we are able to show that the time to escape from extinction is asymptotically equal to the time of the  first death occurring in the limit skeleton. An explicit formula for the density distribution function of the  time to the first death in the limit skeleton is derived.

\section{Infinite lineages}\label{Ssc}
 Consider a sequence of branching processes $\{Z_m (n)\}_{n=0}^\infty$, $m=1,2,\ldots$ with offspring distributions $(p_m(0),p_m(1),\ldots)$ and starting from one particle $Z_m(0)=1$. Assume that the processes are nearly critical with
\begin{align}
&\sum_{k=1}^\infty kp_m(k)=1+\epsilon_m,\;\; \epsilon_m\to0,\;\; m\to\infty,\label{m1}\\
&\sum_{k=2}^\infty k(k-1)p_m(k)\to \sigma^2,\ m\to\infty, \mbox{ for some } \sigma>0,\label{s1}\\
&\sup_m\sum_{k=n}^\infty k^2p_m(k)\to0, \ n\to\infty.\label{ui1}
\end{align}
Condition \eqref{ui1} requires uniform integrability for the sequence of squared offspring numbers and implies that the following equality
\be
\sum_{k=0}^\infty s^kp_m(k)=1-(1+\epsilon_m)(1-s)+(\sigma^2/2-R_m(s))(1-s)^2
\label{R}
\ee
holds with $R_m(s)\to0$ uniformly in $m$ as $s\uparrow 1$.

A natural way of defining a skeleton for branching processes was proposed in \cite{O}, where such processes were considered with $\epsilon_m>0$ for all $m$ (slightly supercritical case). Its survival probability, $Q_m$, according to Lemma 3.3  in \cite{O} satisfies the following well known approximation formula
\be
Q_m\sim2\epsilon_m\sigma^{-2},\ m\to\infty.
\label{bq}
\ee

It is a well known fact that a supercritical branching process can be viewed as a two-type branching process, by distinguishing  among particles with infinite line of descent and particles having finite number of descendants. If we concentrate only in the number of particles with infinite line of descent we arrive at the so-called skeleton process. Conditioning the supercritical process on non-extinction and focussing on infinite lineages we get a new sequence of supercritical branching process $\{X_m(n)\}_{n\ge0}$ with $X_m(0)=1$ which under conditions \eqref{m1}, \eqref{s1},
and \eqref{ui1} is weakly approximated
\be
\{X_m(t/\epsilon_m)\}_{t\ge0}\to \{Y_1(t)\}_{t\ge0},\ m\to\infty
\label{Ya}
\ee
by the Yule process, see Theorem 3.2 in \cite{O}. Recall that the Yule process is a continuous time Markov branching process with particles living exponential times with mean 1 and, at the moment of death, the particles are replaced by two new particles.
The key part of the proof of \eqref{Ya} is to show, using \eqref{R}, that
\be\label{12}
\left\{
\begin{array}{l}
\mathbb P(X_m(1)=1|X_m (0)=1)=1-\epsilon_m+o(\epsilon_m),  \\
 \mathbb P(X_m(1)=2|X_m (0)=1)=\epsilon_m+o(\epsilon_m).
\end{array}
\right.
\ee
Then it remains to check the convergence of the generator of this Markov chain to the generator of the Yule process after the time is scaled accordingly.

The limitation of this definition of a skeleton is that it has no direct extension to the critical or subcritical branching processes. Theorem 2.1 from \cite{O} shows that in the critical case if the branching process is conditioned "on very late extinction" then the limiting skeleton (without any scaling) is a trivial discrete time process $Y(n)\equiv1$, $n=0,1,\ldots$.

In this paper we suggest an alternative approach of defining a skeleton relying on a random marking of the lineages in the family tree. We start by studying in the next section a simple case of exactly critical reproduction.

\section{Critical branching processes with independently marked particles}\label{Scr}

Consider a single type BGW process $Z(n)$ such that its offspring distribution $(p_0,p_1,\ldots)$ has mean $\sum_{k=1}^\infty kp_k=1$ and finite variance $\sigma^2=\sum_{k=2}^\infty k(k-1)p_k$.
Consider the corresponding family tree  and suppose that each vertex in the tree is independently marked with a small probability $\mu_m\to0$, as $m\to\infty$. Any path connecting the root with a marked vertex will be considered as a part of the skeleton. Thus the skeleton is the subtree of the family tree formed by the skeleton paths.
Adapting Proposition 2.1 from \cite{SS} to this case one can show that for a given $\mu_m$ the sequence $\{X_m(n)\}_{n\ge0}$ of numbers of branches in the skeleton forms a BGW process. Next we find the conditional asymptotic structure of the skeleton.

Let again $Q_m=\mathbb P(X_m (0)=1)$ stand for the probability that the skeleton is not empty (at least one particle is marked). Due to the branching property we have
\be
1-Q_m=  (1-\mu_m) \phi(1-Q_m),
\label{Bp}
\ee
where $\phi(s)=\sum_{k=1}^\infty p_ks^k$. Indeed, \eqref{Bp} simply says that the skeleton is empty if and only if the root is not marked and all the daughter subtrees, if any, have empty skeleton.
Using the Taylor expansion of $\phi$ around point $1$ we get
$$ \phi(1-Q_m) =  1- Q_m+ Q_m^2 \sigma^2/2 + o(Q_m^2),$$
and deduce from \eqref{Bp}
\be
Q_m\sim\sigma^{-1} \sqrt{2\mu_m}.
\label{quc}
\ee
Denoting by $\xi_m$ the indicator of the event that the ancestral particle is marked we get
\begin{align*}
\mathbb E\left(r^{\xi_m}s^{X_m (1)};X_m (0)=1\right) = \mathbb E(r^{\xi_m})\mathbb E(s^{X_m (1)})&-\mathbb E(r^{\xi_m}s^{X_m (1)}; X_m(0)=0)\\
\qquad =(r\mu_m+1-\mu_m)\phi(sQ_m+1-Q_m)&-\mathbb P(X_m (0)=0)
\end{align*}
implying that  the offspring distribution of the skeleton particles satisfies
\begin{align*}
\mathbb E\left(r^{\xi_m}s^{X_m (1)}|X_m (0)=1\right)
& = 1- \frac{1-\phi(1-Q_m(1-s))}{Q_m}-(1-r){\mu_m\over Q_m} + o(\sqrt{\mu_m})\\
& = s+ \frac {\sigma^2}{  2} Q_m(1-s)^2 -\sigma\sqrt{2\mu_m}\cdot{1-r\over2} + o(\sqrt{\mu_m})\\
& = (1-\sigma\sqrt{2\mu_m})s + \sigma\sqrt{2\mu_m} \left(\frac 1 2 r+ \frac 1 2 s^2 \right) + o(\sqrt{\mu_m}).
\end{align*}
It follows that using $\tau_m=\sigma\sqrt{2\mu_m}$ we can write
\be\label{12c}
\left\{
\begin{array}{l}
\mathbb P(\xi_m=0, X_m(1)=1|X_m (0)=1)=1-\tau_m+o(\tau_m),  \\
\mathbb P(\xi_m=1, X_m(1)=0|X_m (0)=1)\sim \tau_m/2,  \\
 \mathbb P(\xi_m=0, X_m(1)=2|X_m (0)=1)\sim \tau_m/2.
\end{array}
\right.
\ee

Comparing \eqref{12c} to \eqref{12} we conclude that, if the original branching process produces at least one marked particle, there holds a weak convergence in the Skorokhod sense
\be
\{X_m\big(t/\tau_m\big)\}_{t\ge0}\to \{Y_{0.5}(t)\}_{t\ge0},\quad m\to\infty.
\label{c05}
\ee
Here the limit process is a continuous time Markov branching process with the critical binary splitting:
\begin{itemize}
\item particles live exponential times with parameter 1,
\item at the moment of its death each particle  with probability 0.5 leaves no children and  with probability 0.5 produces two children.
\end{itemize}
 The limit process, $Y_{0.5}(.)$, being a critical branching process will eventually go extinct.

 Relation  \eqref{12c} gives an enhanced interpretation of the limit skeleton \eqref{c05}. All marked particles appearing in the branching process can be associated with the tips of the limit skeleton. In particular the total number of the marked particles, $W_m$, conditioned on $W_m>0$ is asymptotically distributed as the total number of leaves, $W$, in the family tree of the limit skeleton. Due to the branching property we have
 \[W\stackrel{d}{=}1_{\{\nu=0\}}+(W'+W'')\cdot 1_{\{\nu=2\}},\]
 where $\nu$ is the number of offspring of the initial particle in the skeleton, $W'$ and $W''$ are i.i.d. with $W$. In terms of the generating functions we get an equation $\mathbb E\left(s^W\right)=(s+ \mathbb [E\left(s^W\right)]^2)/2$ leading to $\mathbb E\left(s^W\right)=1-\sqrt{1-s}$.

\section{Main result}\label{Snc}
In this section we combine and further develop the two approaches presented in Chapters \ref{Ssc} and \ref{Scr} for a more general model.
Consider a nearly critical reproduction law $\{p_m(k)\}_{k=0}^\infty$ satisfying \eqref{m1}, \eqref{s1}, \eqref{ui1} and allowing for negative $\epsilon_m$. Furthermore, assume that a particle with $k$ offspring is marked with probability  $A_m(k)$. The marking event may depend not only on the number of offspring but also on the whole daughter branching process. For example, the marking rule could be to mark all particles whose total number of descendent exceeds $m$ \cite{S}.
Observe that for the marked near-critical BGW process the total probability for a particle to be marked is given by
\[
\mu_m=\sum_{k=0}^\infty p_m(k)A_m(k).\]
Clearly, the case of Section \ref{Ssc} corresponds to the zero marking probability, $\mu_m\equiv0$,
and the case of  Section \ref{Scr} corresponds to  $\epsilon_m\equiv0$ and $A_m(k)\equiv\mu_m$.

Reconciling the two different definitions of a skeleton given in Sections \ref{Ssc} and \ref{Scr} we next introduce a new definition.

\begin{definition} For a given family tree of a BGW process with marking, the subtree formed by the marked lineages together with infinite lineages will be called the skeleton.
\end{definition}

Clearly, if $\mu_{m}=0$ the skeleton is formed only by the infinite lineages. If $\mu_{m}>0$ any infinite lineage becomes marked and we can think that the skeleton is formed only by the marked lineages.

In this paper we study the asymptotic behavior of the skeleton assuming
\be
\mu_m\to0,\quad m \to \infty,
\label{m3}
\ee
restricting ourselves to the cases when the mean offspring number for the marked particles
\begin{align*}
M_m&=\mu_m^{-1}\sum_{k=1}^\infty kp_m(k)A_m(k)
\end{align*}
satisfies
\be
\limsup_{m \to \infty}M_m<\infty.
\label{Mf}
\ee
By this we exclude such extreme situations as, for example, when $A_m(k)$ is of order $ \mu_mk^2$ for large $k$ and $\sum k^3p_m(k)\to\infty$. Observe also that given  \eqref{ui1} and  \eqref{m3}
\begin{align}
\sum_{k=2}^\infty k(k-1)p_m(k)A_m(k)\to0\label{s3},
\end{align}
which is obtained by using the inequality
\begin{align*}
\sum_{k=2}^\infty k(k-1)p_m(k)A_m(k)\le n^2\mu_m+\sum_{k=n}^\infty k^2p_m(k).
\end{align*}

Let, as before, $Q_{m}$ stand for the probability that the skeleton is not empty. Now we can state our main result claiming that, conditioned on the event that the skeleton is not empty, a weak convergence of the following form holds
\be
\{X_m(t/\tau_m)\}_{t\ge0}\to \{Y_{\lambda}(t)\}_{t\ge0},\quad m\to\infty,
\label{cla}
\ee
for a convenient sequence $(\tau_m)_{m \geq 0}$ and convenient $\lambda \in [0,1]$, generalizing both \eqref{Ya} and \eqref{c05}. Here for a given $\lambda\in[0,1]$ the limit process is a continuous time Markov branching process with binary splitting:
\begin{itemize}
\item particles live exponential times with parameter 1,
\item at the moment of its death each particle  with probability $1-\lambda$ leaves no children and  with probability $\lambda$ produces two children.
\end{itemize}
{\bf Remark.} Importantly, as with  \eqref{c05} in Section \ref{Scr}, by the claiming \eqref{cla} we implicate that asymptotically there is one-to-one correspondence among the marked particles appearing in the branching process and the tips of the limit skeleton $Y_\lambda (.)$. In particular, in  \eqref{Ya} the limit skeleton has no tips implying that under the corresponding time scale we can not expect observing marked particles in the branching process.

\begin{thm}\label{theo}  Under conditions \eqref{m1}, \eqref{s1}, \eqref{ui1},  \eqref{m3},  \eqref{Mf} assuming that there exists a finite or infinite limit
\be
c=\lim_{m\to\infty}\epsilon_m/\sqrt{\mu_m},\label{c}
\ee
\begin{description}
\item[ (i)] if $c=\infty$, then \eqref{bq} and \eqref{Ya} hold,
\item[ (ii)] if $c\in(-\infty,\infty)$,
then
\be
Q_{m}\sim \sqrt{\mu_m}\cdot{c+\sqrt{c^2+2\sigma^2}\over\sigma^2},\label{Qm}
\ee
also  \eqref{cla} holds with $\tau_m= \sqrt{\mu_m}\sqrt{c^2+2\sigma^2}$ and $\lambda={1\over2}+{1\over2}{c\over\sqrt{c^2+2\sigma^2}}$,
\item[ (iii)] if  $c=-\infty$, then
\[Q_{m}\sim {\mu_{m}/|\epsilon_m|},\]
also  \eqref{cla} holds with $\tau_m=|\epsilon_m|$ and $\lambda=0$.
\end{description}
\label{Thm1}
\end{thm}

According to Theorem \ref{theo} there are five different asymptotic regimes for the skeleton of a near-critical BGW process depending on how the deviation from the critical reproduction, $\epsilon_m$, relates to the square root of the marking probability, $\sqrt{\mu_m}$:
\begin{itemize}
\item  in the supercritical case $c=\infty$ with a negligible marking probability the limit skeleton is the Yule process which never dies out,
\item  in the supercritical case $c\in(0,\infty)$ with a balanced marking probability the limit skeleton is a supercritical Markov branching process which dies out with probability $\frac{\sqrt{c^2+ 2\sigma^2}-c}{ \sqrt{c^2+ 2\sigma^2}+c}$  and survives forever with probability $\frac{2c}{ \sqrt{c^2+ 2\sigma^2}+c}$,
\item  if the reproduction law is very close to the purely critical, $c=0$, then the limit skeleton is a critical Markov branching process which dies out with probability one although rather slowly,
\item  in the subcritical case $c\in(-\infty,0)$ with a balanced marking probability the limit skeleton is a subcritical Markov branching process which dies out with probability one,
\item  in the subcritical case $c=-\infty$ with a very small marking probability the limit skeleton is given by a single lineage that dies out after an exponential time.
\end{itemize}

\section{Binomial mutation model}\label{Sbm}

Here we present an important example of a marked branching process based on a decomposable two-type Galton-Watson process modeling a population of individuals with irreversible mutations. The two-type branching process stems from a single wild type individual which produces $k$ offspring with probability $q_m(k)$. Suppose that each daughter of a wild type individual becomes a mutant with probability $\pi_m$ independently of other daughters.

To introduce a marked BGW process we focus only on the wild type individuals and mark those wild type individuals who have at least one mutant daughter. The reproduction law for the marked branching process is given by the distribution for the number of wild type offspring:
\be
p_m(k)= \sum_{l=0}^\infty q_m(k+l){k+l\choose l}(1-\pi_m)^k\pi_m^l,\label{fpm}
\ee
and the conditional marking probabilities $A_m(k)$ are computed using the following relations obtained by splitting \eqref{fpm} in two parts
\begin{align}
p_m(k)(1-A_m(k))&= q_m(k)(1-\pi_m)^k,\label{puk}\\
p_m(k)A_m(k)&= \sum_{l=1}^\infty q_m(k+l){k+l\choose l}(1-\pi_m)^k\pi_m^l.\label{pum}
\end{align}

To ensure that one can use the results from previous section, we need conditions \eqref{m1}-\eqref{ui1} to hold. Therefore we assume that the reproduction law with mutant offspring satisfies
\begin{align}
&\sum_{k=1}^\infty kq_m(k)= 1+\eta_m,\;\; \eta_m \to 0,\label{m4}\\
&\sum_{k=2}^\infty k(k-1)q_m(k) \to \sigma^2,\label{s4}\\
&\sup_m \sum_{k=n}^\infty k^2q_m(k) \to 0, \;\; n \to \infty,\label{ui4}
\end{align}
for some  $\sigma\in(0,\infty)$. Assume also
\be
\pi_m\to0,\quad m\to\infty.
\label{pim}
\ee

\begin{lemma}
Conditions \eqref{m4}, \eqref{s4}, \eqref{ui4},  \eqref{pim}
imply \eqref{m1}, \eqref{s1}, \eqref{ui1},  \eqref{m3},  \eqref{Mf} with
\begin{align}
\mu_m&\sim\pi_m,\quad m\to\infty,\label{mum}\\
M_m&\to\sigma^2,\quad m\to\infty,\label{Mum}.
\end{align}
\end{lemma}
{\bf Remark.} Relation \eqref{Mum} has an interesting implication for our two-type branching process: the mean number of wild type siblings in a family with at least one mutant asymptotically equals the variance of the total offspring number.

\begin{proof}
First observe that   \eqref{pum} entails a useful expression for the marking probability
\begin{align}
\mu_m=\sum_{k=0}^\infty p_m(k)A_m(k)&= \sum_{k=1}^\infty q_m(k)(1-(1-\pi_m)^k).\label{pul}
\end{align}
Clearly, \eqref{pul} and \eqref{m4} yield
\begin{align}
 0\le \mu_m-\pi_m (1+\eta_m)\le \pi_m^2 \sum_{k=1}^\infty k^2q_m(k), \label{mk}
\end{align}
and \eqref{mum} follows from \eqref{m4}, \eqref{ui4} and \eqref{mk}.

Next, due to \eqref{puk} we have
\begin{align*}
0&\le  \sum_{k=2}^\infty k(k-1)q_m(k)-\sum_{k=2}^\infty k(k-1)p_m(k)(1-A_m(k))\\
&=  \sum_{k=2}^\infty k(k-1)q_m(k)(1-(1-\pi_m)^k)\\
&\le n^3\pi_m+\sup_j \sum_{k=n}^\infty k^2q_j(k)
\end{align*}
for any $n\ge2$. Letting here first $m\to\infty$ and then $n\to\infty$, due to \eqref{s4} and \eqref{ui4}, we arrive at
\begin{align*}
&\sum_{k=2}^\infty k(k-1)p_m(k)(1-A_m(k))\to\sigma^2.
\end{align*}
This together with \eqref{s3} implies  \eqref{s1}.
Now, according to \eqref{pum} we have
\begin{align*}
\mu_mM_m\pi_m^{-1}&= \sum_{k=1}^\infty k\sum_{l=1}^\infty q_m(k+l){k+l\choose l}(1-\pi_m)^k\pi_m^{l-1}\\
&= \sum_{j=2}^\infty q_m(j)\sum_{l=1}^j  (j-l){j\choose l}(1-\pi_m)^{j-l}\pi_m^{l-1}\\
&= (1-\pi_m)^{-1}\sum_{j=2}^\infty j (j-1)p_m(j)(1-A_m(j))\\
&\quad+\pi_m\sum_{j=2}^\infty q_m(j)\sum_{l=2}^j  (j-l){j\choose l}(1-\pi_m)^{j-l}\pi_m^{l-2}.
\end{align*}
From here we easily obtain \eqref{Mum}, and therefore   \eqref{Mf}, using   \eqref{s3} and
\begin{align*}
 \sum_{j=2}^\infty q_m(j)\sum_{l=1}^j  (j-l){j\choose l}(1-\pi_m)^{j-l}\pi_m^{l-2}\le  \sum_{j=2}^\infty j(j-1)q_m(j).
\end{align*}

In view of
\[\sum_{k=1}^\infty kp_m(k)(1-A_m(k))=1+\epsilon_m-\mu_mM_m\]
we derive from \eqref{puk} and \eqref{m4}
\begin{align}
0\le \eta_m-\epsilon_m+\mu_mM_m\le  \pi_m \sum_{k=1}^\infty k^2q_m(k).\label{qk}
\end{align}
Combining  \eqref{Mf} and \eqref{qk} we get \eqref{m1}.

To prove \eqref{ui1} we turn to \eqref{fpm} and see that
\begin{align*}
 \sum_{k=n}^\infty k^2p_m(k)&= \sum_{k=n}^\infty k^2\sum_{l=0}^\infty q_m(k+l){k+l\choose l}(1-\pi_m)^k\pi_m^{l}\\
&= \sum_{j=n}^\infty q_m(j)\sum_{l=0}^{j-n}  (j-l)^2{j\choose l}(1-\pi_m)^{j-l}\pi_m^{l}\\
&\le  \sum_{j=n}^\infty j^2q_m(j).
\end{align*}
Thus  \eqref{ui1} is an immediate consequence of  \eqref{ui4}.
\end{proof}

\begin{cor}\label{co}
 For the binomial mutation model satisfying  \eqref{m4}, \eqref{s4}, \eqref{ui4}, \eqref{pim}, and $\eta_m/\sqrt{\pi_m}\to c$ the statements {\bf (i), (ii), (iii)} of Theorem \ref{theo} are valid after $(\epsilon_m,\mu_m)$ are replaced by $(\eta_m,\pi_m)$.
\end{cor}

\section{The sequential mutation model}\label{Ssm}

Our next illustration of Theorem \ref{theo} deals with the sequential mutation model \cite{SS} for a viral population with irreversible mutations which escapes extinction as soon as a target type of viruses is produced. To simplify the discussion we focus mainly on the two-step mutation model, extending the one-step model from Section \ref{Sbm}.

Suppose we have a population of viruses stemming from a single virus which is able to reproduce and mutate giving rise to what we call intermediate type of viruses. The viruses of intermediate type reproduce according to a common law and by mutation generate a new type of viruses which we call the target type. The marking rule for the intermediate type is straightforward: we mark mothers with at least one daughter of the target type.  The wild type marking rule is a bit more complicated: we mark a mother which has at least one  {\it successful mutant daughter} (that is a mutant, of the intermediate type, which has at least one marked descendant in the whole line of descent).

We will assume that the reproduction laws and marking probabilities for both wild type and intermediate type branching processes satisfy conditions of type \eqref{m4}, \eqref{s4}, \eqref{ui4}, \eqref{pim} and are described by triplets $(\eta_m,\sigma^2,\pi_m)$ and $(\hat\eta_m,\sigma_2^2,\hat\pi_m)$ respectively. Furthermore, we assume for simplicity that
$$\hat\pi_m/\pi_m\to\alpha,\ \ m\to\infty,\mbox{ for some }\alpha\in[0,\infty).$$
According to Corollary \ref{co}, if
$$\hat\eta_m\sim c_2\sqrt{\pi_m},\ \ m\to\infty,\mbox{ for some }c_2\in(-\infty,\infty),$$
then
\[\hat Q_{m}\sim \sqrt{\pi_m}\cdot{c_2+\sqrt{c_2^2+2\alpha\sigma_2^2}\over\sigma_2^2},\]
where $\hat Q_m$ is the counterpart of $Q_m$ in Corollary \ref{co} for the skeleton of the intermediate type. Clearly, the intermediate type is supercritical iff $c_2>0$ and $\alpha>0$. Notice that with $\alpha=0$ the skeleton of the intermediate type is the Yule process. The time scale intermediate type skeleton is given by
\begin{align}
 \hat \tau_m&=\sqrt{\pi_m}\sqrt{c_2^2+2\alpha\sigma_2^2}.\label{c2}
\end{align}

Our main interest in the two-step mutation model is of course the limit skeleton leading to the target type.
Therefore, we want to apply  Corollary \ref{co} once again to the branching system with the probability of a successful mutation for the wild type viruses given by
$$\tilde\pi_m=\pi_m\hat Q_m\sim \pi_m^{3/2}\cdot{c_2+\sqrt{c_2^2+2\alpha\sigma_2^2}\over\sigma_2^2}.$$
 Assuming
 $$\eta_m\sim c_1\pi_m^{3/4},\ \ m\to\infty,\mbox{ for some } c_1\in(-\infty,\infty)$$
 we get
$$Q_m\sim \sqrt{\tilde \pi_m}\cdot{c+\sqrt{c^2+2\sigma^2}\over\sigma^2},\quad c={c_1\sigma_2\over\sqrt{c_2+\sqrt{c_2^2+2\alpha\sigma_2^2}}},$$
implying
\be
Q_m\sim \pi_m^{3/4}\cdot{c_1\sigma_2+\sqrt{c_1^2\sigma_2^{2}+2\sigma^2(c_2+\sqrt{c_2^2+2\alpha\sigma_2^2})}\over\sigma^2\sigma_2}.\label{Q2}
\ee
The limit skeleton for the wild type is supercritical iff $c_1>0$.

The time scale for the wild type  type skeletons is given by
\begin{align*}
\tau_m&=\pi_m^{3/4}\sqrt{c_1^2+2\sigma^2\sigma_2^{-2}(c_2+\sqrt{c_2^2+2\alpha\sigma_2^2})},
\end{align*}
which in the considered case is much slower than the time scale of the intermediate type \eqref{c2}. Thus the overall skeleton is given by the wild type skeleton, and the first death in the limit skeleton corresponds to the time of escape from extinction when the first virus of the target type appears.

The above considered case is one of the many possible combination of reproduction-mutation regimes for the two-step mutation model. Without analyzing each of the remaining cases we just point out that there is a situation when both parts of the skeleton live on the same time scale. This is the case when
\begin{align*}
 \hat \pi_m&\sim\alpha \pi_m^{\gamma},\ \ m\to\infty,\mbox{ for some }\gamma\in(1,2),\\
 \hat \eta_m&\sim-\beta \pi_m^{\gamma-1},\ \ m\to\infty,\mbox{ for some }\beta\in(0,\infty),\\
 \eta_m&\sim c_1\pi_m,\ \ m\to\infty,\mbox{ for some } c_1\in(-\infty,\infty).
\end{align*}
Here both time scales are of order $1/\pi_m$.

Turning to the the sequential model with $b-1$ intermediate steps before the target type, we extrapolate the formula \eqref{Q2} to
$Q_m\sim {\rm const}\cdot\pi_m^{1-2^{-b}}$. This prediction should be compared with the strictly subcritical case where one expects $Q_m\sim {\rm const}\cdot\pi_m^{b}$, see Theorem 7.1 in \cite{SS}.

\section{Time to escape}

In the framework of the sequential mutation model of Section \ref{Ssm} it  is crucial to be able to describe the time to escape from extinction.

Let $T_m$ be the time until the first marked particle is observed and put
\[Q(t)=\lim_{m\to\infty}\mathbb P\Big(T_m>{t\over\sqrt{\mu_m(c^2+2\sigma^2)}}\Big).\]
According to Theorem \ref{theo} given that the limit \eqref{c} is finite,  $Q(t)=\mathbb P(T>t)$ is the tail probability  of the time $T$ to the first death in the limit skeleton $Y_\lambda(\cdot)$ with $\lambda={1\over2}+{1\over2}{c\over\sqrt{c^2+2\sigma^2}}$. The branching property of $Y_\lambda(\cdot)$ says that
\be
T=L+\min(T',T'')\cdot 1_{\{\nu=2\}}, \label{TL}
\ee
where $L$ is the exponential life length with mean one, $\nu$ is the number of offspring of the initial particle, $T'$ and $T''$ are i.i.d with $T$.
Due to the branching property \eqref{TL}
\begin{align*}
Q(t)& = P(L>t) + P(T>t, L\leq t)\\
&=e^{- t}+ \lambda\int_0^tQ^2(t-u) e^{- u}du.
\end{align*}
It follows
\begin{align*}
e^{t}Q(t)&=1+ \lambda\int_0^tQ^2(u) e^{u}du.
\end{align*}
Differentiation over $t$ yields a simple differential equation
\begin{align*}
Q'(t)+Q(t)&= \lambda Q^2(t),\quad Q(0)=1
\end{align*}
giving $Q(t)=1/(\lambda+(1-\lambda)e^t)$. Thus
\[\mathbb P\Big(T_m>{t\over\sqrt{\mu_m}}\Big)\to 2\Big(1+{c\over\sqrt{c^2+2\sigma^2}}+\big(1-{c\over\sqrt{c^2+2\sigma^2}}\big)e^{t\sqrt{c^2+2\sigma^2}}\Big)^{-1}\]
and we conclude that the scaled time $T_m\sqrt{\mu_m}$ has the limit density distribution function
\be
\psi(t)= {2(\sqrt{c^2+2\sigma^2}-c)e^{t\sqrt{c^2+2\sigma^2}}\over\Big(1+{c\over\sqrt{c^2+2\sigma^2}}+\big(1-{c\over\sqrt{c^2+2\sigma^2}}\big)e^{t\sqrt{c^2+2\sigma^2}}\Big)^2},\quad t\ge0.\label{den}
\ee
In particular, with $c=0$ we get
\[\psi(t)= 2\sqrt{2}\sigma^2e^{t\sigma\sqrt{2}}\Big(1+e^{t\sigma\sqrt{2}}\Big)^{-2}.\]

\begin{figure}[h!]
\centering
\includegraphics[width=0.8\textwidth]{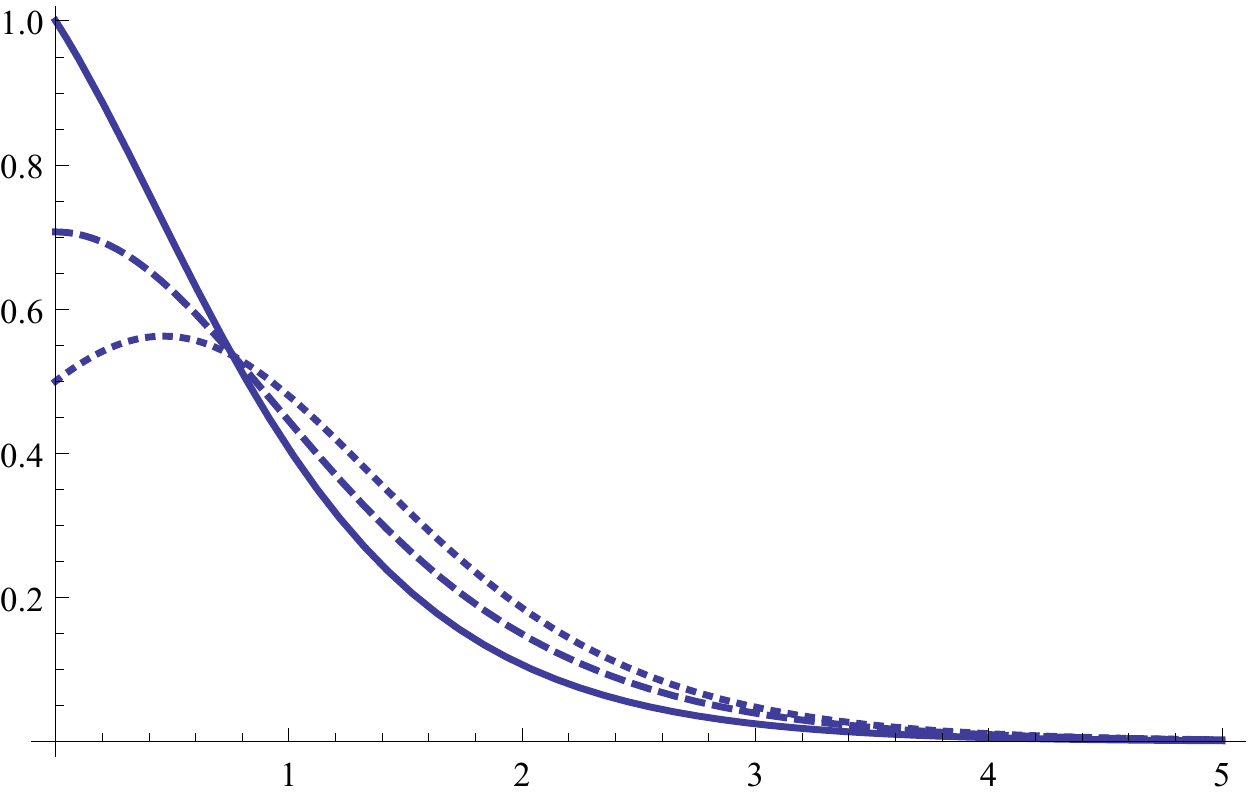}
\caption{
Three density curves given by formula \eqref{den} with $\sigma=1$: subcritical case $c=-0,5$ (solid line), critical case $c=0$ (dashed line), supercritical case $c=0.5$ (dotted line).
}
\label{f1}
\end{figure}

We illustrate the asymptotic density function, $\psi$, by Figure \ref{f1}. In the supercritical case the density curve reaches its a maximum value at
\[t_{\rm max}={1\over\sqrt{c^2+2\sigma^2}} \ln\Big( 1+{2c\over \sqrt{c^2+2\sigma^2}-c}\Big)\]
making the most likely value for the time to escape $T_m$ to be around
\[\hat T_m= \sqrt{{1\over \mu_m(c^2+2\sigma^2)}} \ln\Big( 1+{2c\over \sqrt{c^2+2\sigma^2}-c}\Big).\]

\section{Proof of Theorem \ref{theo}}

Our proof relies on the properties of the probability generating function
$$f_m(r,s)= \sum_{k=0}^\infty p_m(k)s^k \big(1-A_m(k)+ A_m(k)r\big)$$
jointly characterizing the marking status of a particle (through $r$) and its offspring number (through $s$).

\begin{lem} Given conditions \eqref{m1},  \eqref{s1},  \eqref{ui1}, and any sequence $u_m\in(0,1)$ such that $u_m\to 0$, as $m\to \infty$, the following decomposition holds
$$1-f_m (r, 1-u_m)= \mu_m(1-r) + u_m\big(1+\epsilon_m -\mu_m M_m (1-r)\big) -u_m^2\sigma_{r,m} ^2/2,$$
with $\sigma_{r,m} ^2 \to \sigma^2$ as $m\to \infty$ uniformly over $r\in[0,1]$.
\end{lem}
\label{aux_lema}
\begin{proof}
This follows from a Taylor expansion around point $(r,1)$
$$f_m (r, 1-u_m)= f_m(r,1)- u_m \frac{\partial f_m}{\partial s} (r,1)+ \frac{u_m^2}{2} \left( \frac{\partial^2 f_m}{\partial s^2}(r,1)+ R_m(r) \right),$$
where
\begin{align*}
f_m(r,1)& = \sum_{k=0}^\infty p_m(k) (1-A_m(k)+ A_m(k)r)= 1-\mu_m(1-r),\\
\frac{\partial f_m}{\partial s} (r,1) &= \sum_{k=1}^\infty k p_m(k) (1-A_m(k)+ A_m(k)r)= 1+ \epsilon_m -\mu_m M_m(1-r),\\
\frac{\partial^2 f_m}{\partial s^2} (r,1) &= \sum_{k=2}^\infty k (k-1)p_m(k) (1-A_m(k)+ A_m(k)r),
\end{align*}
and
\[
R_m(r)= \sum_{k=2}^\infty k (k-1)p_m(k) (1-A_m(k)+ A_m(k)r) (1-\theta_m^{k-2}),\,
\]
for some $\theta \in (1-u_m, 1)$. Indeed, since for any $n\ge3$
\begin{align*}
R_m(r)&\le \sum_{k=2}^\infty k^2p_m(k)(1-(1-u_m)^{k-2})\\
&\le u_mn^{3}+\sum_{k=n}^\infty k^2p_m(k),
\end{align*}
condition \eqref{ui1} implies $R_m(r)\to0$.
It remains to apply \eqref{s1} and \eqref{s3}.
\end{proof}

The skeleton is empty if the initial particle is not marked and all her children produce empty skeletons
$$1-Q_m = f_m(0, 1-Q_m).$$
Using Lemma 8.1, with $r=0$ and $u_m= Q_m$, we obtain a quadratic equation
$$ \sigma_{0,m}^2 Q_m^2 -2 Q_m (\epsilon_m -\mu_m M_m) -2\mu_m= 0,$$
entailing
\be \label{2-sol}
Q_m = \frac{\epsilon_m - \mu_m M_m + \sqrt {(\epsilon_m -\mu_mM_m)^2 + 2\sigma_{0,m}^2\mu_m}}{\sigma_{0,m}^2},
\ee
where $\sigma_{0,m} \to \sigma$. Once again applying Lemma 8.1  now to the right hand side of
\begin{eqnarray*}
\mathbb E[r^{\xi_m} s^{X_m(1)}; X_m(0)=1] & = &   E[r^{\xi_m} s^{X_m(1)}]- P[X_m(0)=0]\\
& =  & f_m(r, s Q_m + 1- Q_m) - (1-Q_m)
\end{eqnarray*}
with $u_m= 1-Q_m(1-s)$, we get
\begin{align}\label{asymp_skelet}
\mathbb E[r^{\xi_m} s^{X_m(1)}|X_m(0)=1] & =  \left( 1+\epsilon_m -Q_m \sigma^2 \right)s +Q_m s^2\sigma^2/2 + \mu_mQ_m^{-1}r
 \\\nonumber
&\qquad-\mu_mQ_m^{-1} -\epsilon_m  + Q_m \sigma^2/2 + o(Q_m)+O(\mu_m ).
\end{align}

Now we are ready to verify the statements of Theorem $\ref{theo}$ case by case.

{\bf Case (i)} If $c=\infty$, then  is $\mu_m=o(\epsilon_m)$ and  $\eqref{2-sol}$ yields \eqref{bq}.
Furthermore, $\mu_mQ_m^{-1} = o(\epsilon_m)$ and  $\eqref{asymp_skelet}$ gives
\be
\mathbb E[r^{\xi_m} s^{X_m(1)}|X_m(0)=1] = (1-\epsilon_m)s+ \epsilon_m s^2
\ee
as desired.

{\bf Case (ii)}  If $\epsilon_m \sim c\sqrt{\mu_m}$ for  $c\in (-\infty, +\infty)$,  then $\eqref{2-sol}$ yields \eqref{Qm} as we can neglect the terms involving $\mu_mM_m$.
Furthermore, $\mu_mQ_m^{-1} \sim \frac{\sqrt{\mu_m}\sigma^2}{c+ \sqrt{c+ c\sigma^2}}$ and the last terms in $\eqref{asymp_skelet}$ are negligible
\begin{align*}
-\mu_mQ_m^{-1}& -\epsilon_m + Q_m \sigma^2/2 + o(Q_m)+O(\mu_m )\\
&=   -\frac{ \sqrt{\mu_m}\sigma^2}{c+ \sqrt{c+ c\sigma^2}}\, - c\sqrt{\mu_m} \,+  \sqrt{\mu_m} \frac{c+ \sqrt{c^2+ 2\sigma^2}}{2} + o(\sqrt{\mu_m}) =o(\sqrt{\mu_m}).
\end{align*}
We conclude
\begin{align*}
&\mathbb E[r^{\xi_m} s^{X_m(1)}|X_m(0)=1] =\left (1- \sqrt{\mu_m} \sqrt{c^2 + 2\sigma^2}\right) s \\
&\quad+  \sqrt{\mu_m} \sqrt{c^2 + 2\sigma^2} \left[ \left(\frac{1}{2} - \frac{1}{2} \frac{c}{\sqrt{c^2 + 2\sigma^2}}\right) r +  \left(\frac{1}{2} + \frac{1}{2} \frac{c}{\sqrt{c^2 + 2\sigma^2}}\right)s \right]+ o(\sqrt{\mu_m})
\end{align*}
as desired.

{\bf Case (iii)}  If $c= -\infty$, then $a_m:=\epsilon_m-\mu_m M_m$ takes negative values for sufficiently large $m$ and we derive from \eqref{2-sol}
\[
Q_m = \frac{2\sigma_{0,m}^2 \mu_m}{\sigma_{0,m}^2\left(\sqrt{a_m^2+2\sigma_{0,m}^2 \mu_m}-a_m\right)}\sim \frac{\mu_m}{ |\epsilon_m|}.
\]
as stated. It remains to see that \eqref{asymp_skelet} gives
$$\mathbb E[r^{\xi_m} s^{X_m(1)}| X_m(0)=1] = (1-|\epsilon_m|)s + |\epsilon_m|r+ o(\epsilon_m),$$
since $\mu_mQ_m^{-1} \sim |\epsilon_m|$, $\mu_m= o(|\epsilon_m|)$ and $Q_m=o(|\epsilon_m|)$.

\acks
The research of Serik Sagitov was supported by the Bank of Sweden Tercentenary Foundation and the Swedish Research Council grant 621-2010-5623. The research of Maria Concei\c{c}\~{a}o Serra was financed by FEDER Funds through "Programa
Operacional Factores de Competitividade - COMPETE" and by
Portuguese Funds through FCT -  "Funda\c{c}\~{a}o para a Ci\^{e}ncia e a Tecnologia",
within the Project Est-C/MAT/UI0013/2011.

\end{document}